\documentclass[reqno,a4paper]{amsart}

\usepackage[utf8]{inputenc}
\usepackage[T1]{fontenc}

\usepackage[pagebackref]{hyperref}
\usepackage{amsfonts}
\usepackage{bbm}
\usepackage{mathtools}
\usepackage{enumerate}
\usepackage[left=2.5cm,right=2.5cm]{geometry}
\allowdisplaybreaks

\newtheorem{theorem}{Theorem}[section]

\newtheorem{cor}[theorem]{Corollary} 
\newtheorem{prop}[theorem]{Proposition} 

\theoremstyle{definition}
\newtheorem{defi}[theorem]{Definition}

\theoremstyle{remark}
\newtheorem{example}[theorem]{Example} 
\newtheorem{remark}[theorem]{Remark}

\numberwithin{equation}{section}

\newcommand{\ub}[1]{^{(#1)}}
\def\ii{{\bf i}}

\newcommand{\ind}[1]{\mathbbm{1}_{#1}} 
\newcommand{\E}{\mathbb{E}} 
\newcommand{\p}{\mathbb{P}} 
\newcommand{\R}{\mathbb{R}} 

\DeclareMathOperator{\Ent}{Ent}	
\DeclareMathOperator{\Med}{Med}	
\DeclareMathOperator{\Var}{Var}	

\newcommand*{\cTbar}{\overline{\mathcal{T}}}
\newcommand*{\Tbar}{\overline{\mathbf{T}}\vphantom{\mathbf{T}}}

\author{Rados{\l}aw Adamczak} %
\address[RA]{Institute of Mathematics of the Polish Academy of Sciences \\
\& Institute of Mathematics, University of Warsaw } %
\email{R.Adamczak@mimuw.edu.pl}

\author{Micha{\l} Kotowski}
\address[MK]{Institute of Mathematics, University of Warsaw } %
\email{Michal.Kotowski1@gmail.com}

\author{Bart{\l}omiej Polaczyk}
\address[BP]{Faculty of Mathematics, Informatics and Mechanics, University of Warsaw}
\email{B.Polaczyk@student.uw.edu.pl}

\author{Micha{\l} Strzelecki}
\address[MS]{Institute of Mathematics, University of Warsaw } %
\email{M.Strzelecki@mimuw.edu.pl}

\thanks{Research partially supported by the National Science Centre, Poland, grant
no. 2015/18/E/ST1/00214 (RA, MK), 2015/19/N/ST1/00891 (MS) and 2017/24/T/ST1/00323 (doctoral scholarship of MS)}

\title{A note on concentration for polynomials in the Ising model}

\begin{document}

\begin{abstract}
We present precise multilevel exponential concentration inequalities for polynomials in Ising models satisfying the Dobrushin condition. The estimates have the same form as two-sided tail estimates for polynomials in Gaussian variables due to Lata{\l}a. In particular, for quadratic forms we obtain a Hanson-Wright type inequality.

We also prove concentration results for convex functions and estimates for nonnegative definite quadratic forms, analogous as for quadratic forms in i.i.d.\ Rademacher variables, for more general random vectors satisfying the approximate tensorization property for entropy.

\end{abstract}

\maketitle

\section{Introduction}

Since its introduction in \cite{1925ZPhy...31..253I}, the Ising model has been a source of numerous mathematical questions. In addition to its physical importance it is appealing to mathematicians, providing an easy to formulate, yet challenging model of dependent random variables and serving as testing ground for many probabilistic ideas. Recently in the context of finite graphs, the Ising model attracted also attention of statisticians and theoretical computer scientists interested e.g, in estimating the parameters of the model, learning the underlying graph structure or testing some properties of the model in a computationally efficient way (see e.g., \cite{MR3775918,MR2943079,2014arXiv1411.1434S,2018arXiv180606887D}). In particular in the last two decades several authors studied the Ising model from the point of view of concentration of measure phenomena, see e.g., \cite{MR2002993,MR2288072,MR2508781,MR3711609}. While most effort has been devoted to concentration inequalities for functions satisfying appropriate Lipschitz type conditions, recently several papers appeared related to variance bounds or stronger, exponential concentration inequalities for polynomials, to mention the work by Daskalakis, Dikkala, Kamath \cite{MR3775918, 2017arXiv171004170D}, Gheissari, Lubetzky and Peres \cite{2017arXiv170600121G}, G\"otze, Sambale and Sinulis \cite{2018arXiv180106348G}. Motivation for these developments ranged from statistical and algorithmic (efficient discrimination between samples drawn from an Ising model and i.i.d.\ samples) to purely probabilistic ones (searching for counterparts of inequalities known in the i.i.d.\ case).

In this note we complement the results proved in the aforementioned papers, by obtaining exponential inequalities for polynomials of the same form as in bounds for polynomials in independent Gaussian (or more generally subgaussian) random variables, which were introduced originally by Lata{\l}a \cite{MR2294983} and subsequently studied e.g., by Adamczak and Wolff \cite{MR3383337}. Such inequalities are expressed in terms of appropriate injective tensor product norms of averaged derivatives of the polynomials in question and in the Gaussian case are known to be optimal up to constants depending only on the degree of the polynomial. Optimality is understood here in a strong sense -- the inequalities can be up to constants reversed. Moreover they are known to imply other, more classical inequalities for multilinear forms, such as Bonami-Nelson inequalities \cite{MR0283496,MR0343815}. In particular, for polynomials of degree $d$ they provide multilevel type concentration, of the form $\exp(-ct^2)$ for small values of $t$ up to $\exp(-c't^{2/d})$ for larger values (as opposed to the inequalities from the aforementioned results for the Ising model, which do not yield precise multilevel concentration but rather give weaker bounds of the form $\exp(-c''t^{2/d})$ for all $t$). As a consequence our estimates imply the previous ones and provide a more accurate description of the tail behavior.

Our approach is similar to the one by G\"otze, Sambale and Sinulis in that it builds on general Aida--Stroock type moment estimates those Authors obtained for the Ising model in \cite{2018arXiv180106348G}, and uses them as a tool in an inductive argument. However, the details are different, while in \cite{2018arXiv180106348G} one works with moments of Euclidean norms of discrete iterated derivatives of multilinear forms, we adapt an argument from \cite{MR3383337}, linearizing the Euclidean norms with an auxiliary Gaussian sequence, which allows us to treat general functions (seen by the Fourier--Walsh theory as tetrahedral polynomials) and pass from discrete gradients to classical derivatives.

The argument we present may be seen as a method of reduction of concentration properties for polynomials from the Ising model to the i.i.d.\ Gaussian case. Since the random variables considered in the Ising model take only values $\pm1$, one could expect a similar reduction to polynomials in i.i.d.\ Rademacher variables. We are able to obtain such estimates for positive definite quadratic forms, by passing through concentration properties for convex functions which are of independent interest.

\medskip

The organization of the article is as follows. First, in Section \ref{sec:main-results} we present our main result (Theorem \ref{thm:polynomial-Ising}) and discuss its relation with known inequalities for the Ising model as well as with the estimates for the i.i.d.\ case. Next, in Section \ref{sec:tensorization} we discuss the approximate tensorization of entropy (as studied recently by Marton \cite{2015arXiv150702803M} and Caputo, Menz, Tetali \cite{MR3434252}) and Aida--Stroock type moment estimates obtained  by G\"otze, Sambale and Sinulis \cite{2018arXiv180106348G}. Using these tools, in Section \ref{sec:proof-Ising} we present the proof of the main result. The final Section \ref{sec:improved-concentration} presents estimates for convex functions and Rademacher-type inequalities for quadratic forms.

\section{Gaussian type inequality for polynomials}\label{sec:main-results}

\subsection{Basic definitions and notation}
Let us begin by introducing the general form of the Ising model on a finite set.

\begin{defi}[Ising model]
Let $n$ be a positive integer and let $\mu$ be the measure on $\{-1,1\}^n$, having density with respect to the uniform distribution of the form
\begin{align}\label{eq:Ising}
\mu(\sigma) = \frac{1}{Z}\exp\Big(\frac{1}{2}\sum_{i,j=1}^n J_{ij}\sigma_i\sigma_j - \sum_{i=1}^n h_i\sigma_i\Big),
\end{align}
for any $\sigma \in \{-1,1\}^n$, where $J = (J_{ij})_{i,j\le n}$ is a symmetric matrix with vanishing diagonal, $h = (h_i)_{i\le n} \in \R^n$ and $Z$ is a normalizing constant.
\end{defi}

In physical terms the coupling matrix $J$ corresponds to interactions between particles and the vector $h$ describes an external field. The order of magnitude of the constants $J_{ij}$ reflects the temperature (the higher the temperature the smaller the coefficients, which corresponds to weaker interactions), however as our results will be expressed solely in terms of the coefficients $J_{ij}$ and $h_i$, we will not incorporate the temperature into the notation.

\medskip

To obtain concentration inequalities, one needs some control over the coupling constants  and the external field, which will allow for sufficient proximity to the i.i.d.\ case. The conditions we will impose on the model are classical and in the context of concentration of measure appeared already in \cite{2015arXiv150702803M, 2017arXiv170600121G,2018arXiv180106348G}.

\medskip

\paragraph{\bf Main assumptions}
We will assume that
\begin{align}\label{eq:h-bound}
\max_{i\le n} |h_i|\le \alpha
\end{align}
 and the coupling constants satisfy Dobrushin's condition
\begin{align}\label{eq:Dobrushin}
  \max_{i\le n}\sum_{j=1}^n |J_{ij}| \le 1 - \rho
\end{align}
for some $\rho > 0$.
\medskip

In order to formulate concentration of measure estimates for polynomials of the Ising model, which correspond to inequalities obtained by Lata{\l}a for polynomials in independent Gaussian random variables, we will need to introduce a family of injective tensor product norms on $d$-index matrices ($d$-tensors).

To provide transparent notation for multi-indices we will use the following convention. For a positive integer $n$ we will denote $[n] = \{1,\ldots,n\}$. The cardinality of a set $I$ will be denoted by $|I|$. For $\ii = (i_1,\ldots,i_d) \in [n]^d$ and $I\subseteq[d]$ we write $\ii_{I}=(i_{k})_{k\in I}$. We will also denote $|\ii| = \max_{j\le d} {i_j}$ and $|\ii_I| = \max_{j \in I} i_j$. We will often deal with homogeneous polynomials, defined in terms of multi-indexed matrices (tensors). We will say that a $d$-indexed matrix $A = (a_\ii)_{\ii \in [n]^d}$ is symmetric if for every permutation $\sigma$ of the set $[d]$ and every $\ii = (i_1,\ldots,i_d) \in [n]^d$, we have $a_\ii = a_{i_{\sigma(1)},\ldots,i_{\sigma(d)}}$. When $d$ is fixed, we will write simply $A = (a_\ii)_{|\ii|\le n}$. We will say that a $d$-indexed matrix $A = (a_\ii)_{|\ii|\le n}$ has vanishing generalized diagonals if $a_\ii = 0$ for all $\ii = (i_1,\ldots,i_d)$ such that there exist $k\neq l$ with $i_k = i_l$.

\medskip

Let now $P_d$ be the set of partitions of $[d]$ into nonempty, pairwise disjoint sets. For a partition $\mathcal{I} =\{I_1,\ldots,I_k\} \in P_d$, and a $d$-indexed matrix $A = (a_\ii)_{\ii \in [n]^d}$, define
\begin{align}\label{eq:Gaussian_norm_def}
\|A\|_{\mathcal{I}}=\sup\Big\{\sum_{\ii\in [n]^d} a_{\ii}\prod_{l=1}^k x\ub{l}_{\ii_{I_l}}\colon
\|(x\ub{l}_{\ii_{I_l}})\|_2\leq 1, 1\leq l\leq k \Big\},
\end{align}

where $\|(x_{\ii_{I_l}})\|_2 = \sqrt{\sum_{|\ii_{I_l}|\le n} x_{\ii_{I_l}}^2}$. Thus, e.g.,
\begin{align*}
\|(a_{ij})_{i,j\le n}\|_{\{1,2\}}&= \sup\Big\{ \sum_{i,j\le n} a_{ij}x_{ij}\colon \sum_{i,j\le n} x_{ij}^2 \le 1\Big\} = \sqrt{\sum_{i,j\le n}a_{ij}^2} = \|(a_{ij})_{i,j\le n}\|_{HS},\\
\|(a_{ij})_{i,j\le n}\|_{\{1\}\{2\}}&= \sup\Big\{ \sum_{i,j\le n} a_{ij}x_iy_j\colon \sum_{i\le n} x_{i}^2\le 1,\sum_{j\le n}y_j^2 \le 1\Big\} = \|(a_{ij})_{i,j\le n}\|_{\ell_2^n\to \ell_2^n},\\
\|(a_{ijk})_{i,j,k\le n}\|_{\{1,2\} \{3\}} &= \sup\Big\{ \sum_{i,j,k\le n} a_{ijk}x_{ij}y_k\colon \sum_{i,j\le n} x_{ij}^2\le 1,\sum_{k\le n}y_k^2 \le 1\Big\}.
\end{align*}

Note that for simplicity in the notation we skip the outer brackets in the subscript and write e.g., $\|\cdot\|_{\{1\}\{2\}}$ instead of $\|\cdot\|_{\{\{1\}\{2\}\}}$.

In particular for $d=2$, $\|\cdot\|_{\{1,2\}}$ and $\|\cdot\|_{\{1\}\{2\}}$ coincide with the Hilbert--Schmidt and operator norm of a matrix respectively. Note that for every $d$ and $\mathcal{I} \in P_d$ we have $\|A\|_{\mathcal{I}} \le \|A\|_{\{[d]\}} = \sqrt{\sum_{|\ii|\le n} a_\ii^2}$. The norm $\|A\|_{\{[d]\}}$ can be considered a counterpart of the Hilbert--Schmidt norm for higher order tensors.

\medskip

We will use the standard notation $\|X\|_p = (\E |X|^p)^{1/p}$ for a random variable $X$. Sometimes, when dealing with independent random variables $X,Y$, we will use the notation $\E_X$ for the expectation with respect to the variable $X$ (i.e., the conditional expectation given $Y$).
\medskip

In what follows, we will write e.g., $c_a, C_a$ or $c_{a}(b)$ to denote constants depending only on the parameters $a$ or $a,b$ respectively. The values of such constants may change between occurrences.

\medskip
By the Fourier--Walsh expansion (see e.g., \cite{MR3443800}), every function $f\colon \{-1,1\}^n \to \R$ can be written in a unique way as a tetrahedral polynomial, i.e., a polynomial which is affine with respect to every variable (in particular the degree of the polynomial is at most $n$). Therefore in what follows we will restrict our attention to this representation. In particular, when we speak about gradients $\nabla f$ or higher order derivatives $\nabla^k f$, we always think of the usual derivatives of the polynomial function on $\R^n$ given by the tetrahedral representation of $f$.

\subsection{Main result}

The main result of this section is
\begin{theorem}\label{thm:polynomial-Ising} Let $\mu$ be defined by \eqref{eq:Ising} and assume that $|h|_\infty \le \alpha$ and the Dobrushin condition \eqref{eq:Dobrushin} holds. Let $X$ be a random vector distributed according to $\mu$. Then for $d \ge 1$ there exist constants $c_d = c_d(\alpha,\rho)$, such that for any tetrahedral polynomial $f\colon \{-1,1\}^n\to \R$ of degree $d$, and any $t > 0$,
\begin{align}\label{eq:polynomial-Ising}
  \p(|f(X) - \E f(X)| \ge t) \le 2\exp\Big(-c_d\min_{1\le k \le d} \min_{\mathcal{I} \in P_k} \Big(\frac{t}{\|\E \nabla^k f(X)\|_\mathcal{I}}\Big)^{2/|\mathcal{I}|}\Big).
\end{align}
\end{theorem}

\begin{remark}
If $h= 0$ and $f(x) = \sum_{|\ii|\le n} a_\ii \prod_{l=1}^d x_{i_l}$ for some symmetric $d$-indexed matrix with vanishing generalized diagonals, then for any $\mathcal{I}\in P_k$ we can estimate
\begin{align*}
  \|\E \nabla^k f(X)\|_\mathcal{I}^2 & \le \Big(\frac{d!}{(d-k)!}\Big)^2\sum_{i_1,\ldots,i_{k}=1}^n \Big(\sum_{i_{k+1},\ldots,i_d=1}^n a_\ii \E \prod_{l=k+1}^d X_{i_l}\Big)^2\\
  &\le \Big(\frac{d!}{(d-k)!}\Big)^2\sum_{i_1,\ldots,i_{k}=1}^n \E \Big(\sum_{i_{k+1},\ldots,i_d=1}^n a_\ii \prod_{l=k+1}^d X_{i_l}\Big)^2.
\end{align*}
By (the proof of) Lemma 3.1. in \cite{2017arXiv170600121G}, $\E (\sum_{i_{k+1},\ldots,i_d=1}^n a_\ii \prod_{l=k+1}^d X_i)^2 \le C_d \max_\ii |a_\ii|^2 n^{d-k}$, as a consequence $\|\E \nabla^k f(X)\|_\mathcal{I} \le \widetilde{C}_d n^{d/2}\max_{\ii} |a_\ii|$ and
Theorem \ref{thm:polynomial-Ising} (after adjustment of constants) implies
\begin{align}\label{eq:l-infty-bound}
  \p(|f(X) - \E f(X)| \ge t) \le 2\exp\Big(-\widetilde{c}_d \frac{t^{2/d}}{n\max_\ii |a_\ii|^{2/d}}\Big)
\end{align}
for some new constant $\widetilde{c}_d = \widetilde{c}_d(\alpha,\rho)$. This inequality was proved by G\"otze, Sambale, Sinulis in \cite{2018arXiv180106348G} and earlier, up to some additional logarithmic in $n$ factors in the exponent by Gheissari, Lubetzky and Peres in \cite{2017arXiv170600121G}.
\end{remark}
\begin{remark}
An  inequality analogous to \eqref{eq:polynomial-Ising} for polynomials in Gaussian variables has been obtained in \cite{MR3383337}. In this case the inequality can be reversed up to numerical constants depending only on $d$ (in front of and inside the exponent). The proof relied on a reduction to the special case of tetrahedral multilinear forms in independent standard Gaussian variables obtained by Lata{\l}a \cite{MR2294983}. The proof of Theorem \ref{thm:polynomial-Ising} presented below is a simple adaptation of this idea. The inequality \eqref{eq:polynomial-Ising} is also known to hold for polynomials in i.i.d.\ subgaussian random variables \cite{MR3383337}.
\end{remark}

\begin{example}
In this example and the following ones we will let $c$ denote a constant which may depend on the parameters $\alpha$ and $\rho$. Its value may change between occurrences.

For $f(x) = \sum_{i,j=1}^n a_{ij} x_ix_j$ we obtain
\begin{displaymath}
  \p(|f(X) - \E f(X)| \ge t) \le 2\exp\Big(-c\min\Big(\frac{t^2}{\|A\|_{HS}^2 + \sum_{i=1}^n (\sum_{j=1}^n a_{ij}\E X_j)^2},\frac{t}{\|A\|_{\ell_2^n\to\ell_2^n}}\Big)\Big).
\end{displaymath}

Note that if $h = 0$ then $\E X_i = 0$ (as the distribution of $X$ is symmetric) and the right hand side simplifies to
\begin{align}\label{eq:Hanson-Wright-123}
  \p(|f(X) - \E f(X)| \ge t) \le \le 2\exp\Big(-c\min\Big(\frac{t^2}{\|A\|_{HS}^2},\frac{t}{\|A\|_{\ell_2^n\to\ell_2^n}}\Big)\Big),
\end{align}
   which gives a counterpart of the Hanson-Wright inequality known for quadratic forms in independent sub-Gaussian random variables \cite{MR0279864}, which turned out to be useful e.g., in random-matrix theory and statistics (see e.g., \cite{vershynin_2018}). A version of this inequality for strongly mixing Ising models on a lattice was proved by Marton \cite{MR2002993}.

If one further estimates the operator norm by the Hilbert--Schmidt norm, one obtains a weaker tail bound of the form
\begin{align}\label{eq:Beckner-Bonami}
\p(|f(X) - \E f(X)| \ge t)  \le 2\exp\Big(-c\frac{t}{\|A\|_{HS}}\Big).
\end{align}
For quadratic form in independent Rademacher variables such an inequality (together with counterparts for higher order forms) were for the first time established by Bonami \cite{MR0283496}, Beckner \cite{MR0385456} and Gross \cite{MR0420249} in the context of hypercontractivity of semigroups. A counterpart of \eqref{eq:Beckner-Bonami} for quadratic forms of the Ising model has been recently obtained in \cite{2018arXiv180106348G}.
\end{example}

\begin{example}\label{ex:degree-3}

Consider now
$f(x) = \sum_{1\le i,j,k \le n} a_{ijk}x_ix_jx_k$, where $A = (a_{ijk})_{i,j,k\le n}$ is symmetric with vanishing generalized diagonals. Assume also that $h=0$.

In this case Theorem \ref{thm:polynomial-Ising} gives
\begin{align}\label{eq:degree-3-example}
  &\p(|f(X)| \ge t) \\
  &\le 2\exp\Big(-c\min\Big(\frac{t^2}{\|A\|_{\{1,2,3\}}^2+\sum_i (\sum_{jk} a_{ijk}\E X_j X_k)^2},\frac{t}{\|A\|_{\{1,2\}\{3\}}},\frac{t^{2/3}}{\|A\|^{2/3}_{\{1\}\{2\}\{3\}}} \Big)\Big),
\end{align}
where we again used the equality $\E X_i = 0$.

One can wonder if it is possible to obtain estimates just in terms of the norm $\|A\|_{\{1,2,3\}}$, e.g., of the form $\p(|f(X)| \ge t) \le 2 \exp(-ct^{2/3}/\|A\|_{\{1,2,3\}}^{2/3})$ as in the independent case or in the case of quadratic forms discussed above.
Clearly, this is true if one can estimate the quantity $\sum_{i,j\le n} (\E X_iX_j)^2$ by a constant independent of $n$. However it turns out that if one assumes only the Dobrushin condition \eqref{eq:Dobrushin}, then the coefficient $\sum_i (\sum_{jk} a_{ijk}\E X_i X_j)^2$ in general cannot be discarded. To see this consider the Ising model on the one-dimensional interval, e.g., with $J_{i,i+1} = J_{i+1,i} = 1/3$
for $i=1,\ldots,n-1$ and $J_{ij}=0$ otherwise. In this case \eqref{eq:Dobrushin} is clearly satisfied with $\rho = 1/3$ and the Hamiltonian is of the form $-\frac{1}{3}\sum_{i=1}^{n-1}\sigma_{i}\sigma_{i+1}$. Since under the uniform measure on the discrete cube, $\sigma_1$ and the products $\sigma_i\sigma_{i+1}$, $i=1,\ldots,n-1$ are independent Rademacher variables, one can see that under the measure $\mu$ given by \eqref{eq:Ising}, the products $\sigma_i\sigma_{i+1}$, are i.i.d.\ random variables with distribution
\begin{align*}
  \mu(\sigma_i\sigma_{i+1} = 1) &= 1 - \mu(\sigma_i\sigma_{i+1} = -1) = \frac{1}{1 + e^{-2/3}}.
\end{align*}
In particular $\E_\mu \sigma_{i}\sigma_{i+1} = a>0$ is independent of $n$. Consider now a symmetric $3$-indexed matrix $A = (a_{ijk})_{i,j,k \le n}$ with vanishing generalized diagonals and coefficients $a_{ijk}$ defined for $i < j < k$ with the formula
\begin{displaymath}
a_{ijk} = \ind{\{i<j=k-1\}}.
\end{displaymath}
One can easily see that $\|A\|_{\{1,2,3\}}$ is of order $n$ (as $A$ has $O(n^2)$ nonzero coefficients). However, if $X$ is distributed according to $\mu$, then for $f(X) = \sum_{i,j,k\le n}a_{ijk}X_iX_jX_k$, $\sqrt{\Var f(X)} = \|f(X)\|_2$ is of order $n^{3/2}$ (as can be checked by using the equality $\E X_iX_jX_k = 0$, expanding the product $\E f(X)^2$ and performing some elementary combinatorics). This shows that an estimate of the form $\p(|f(X)|\ge t) \le 2\exp(-c (t/\|A\|_{\{1,2,3\}})^\kappa)$ cannot hold with $c,\kappa$ independent of $n$. Of course, by \eqref{eq:l-infty-bound} we do have the inequality $\p(|f(X)| \ge t) \le 2\exp(-ct^{2/3}/n)$. However, \eqref{eq:degree-3-example} leads to an improvement of this inequality. As one can easily check $\|A\|_{\{1,2\},\{3\}}$ is of the same order as $\|A\|_{\{1,2,3\}}$, i.e., of order $n$, $\|A\|_{\{1\}\{2\}\{3\}}$ is of order $\sqrt{n}$ and $\sum_i (\sum_{jk} a_{ijk}\E X_j X_k)^2$ is of order $n^{3}$. Together with some elementary calculations, this gives
\begin{displaymath}
  \p(|f(X)| \ge t)\le 2\exp\Big(-c\min\Big(\frac{t^2}{n^3}, \frac{t^{2/3}}{n^{1/3}}\Big)\Big).
\end{displaymath}
\end{example}

\begin{remark}
In \cite{2018arXiv180106348G}, G\"otze, Sambale and Sinulis propose an appropriate re-centering of a $d$-linear form $f(x) = \sum_{|\ii|\le n} a_\ii \prod_{k=1}^d x_{i_k}$ by a polynomial of lower degree, resulting in a (nonhomogeneous) polynomial $f_{A,d}$ such that $\E \nabla^i f(X) = 0$ for all $i=1,\ldots,d-1$. They prove that then
\begin{displaymath}
  \p(|f_{A,d}(X) - \E f_{A,d}(X)| \ge t) \le 2\exp\Big(-c\frac{t^{2/d}}{\|A\|_{\{1,\ldots,d\}}^{2/d}}\Big)
\end{displaymath}
with some $c = c(\alpha,\rho,d)$.
Since Theorem \ref{thm:polynomial-Ising} applies to general (not necessarily homogeneous) polynomials, it yields a refinement of the above inequality, of the form
\begin{displaymath}
  \p(|f_{A,d}(X) - \E f_{A,d}(X)| \ge t) \le 2\exp\Big(-c\min_{\mathcal{I}\in P_d} \Big(\frac{t}{\|A\|_{\mathcal{I}}}\Big)^{2/|\mathcal{I}|}\Big).
\end{displaymath}

Example \ref{ex:degree-3} shows that in general one cannot eliminate passing to $f_{A,d}$, i.e., the above inequality may not hold for the original polynomial $f$.
\end{remark}

\section{Approximate tensorization of entropy and moment estimates}\label{sec:tensorization}

In this section we will present basic tools (coming mostly from the work by Marton \cite{2015arXiv150702803M} and G\"otze, Sambale, Sinulis \cite{2018arXiv180106348G}) which we will need for the proof of Theorem \ref{thm:polynomial-Ising} and also in Section \ref{sec:improved-concentration} to obtain concentration for convex functions and Rademacher-type bounds for quadratic forms.

Let $\mathcal{X} = \prod_{i=1}^n \mathcal{X}_i$, where $\mathcal{X}_i$ are Polish spaces with their Borel $\sigma$-fields and let $\mu$ be a probability distribution on $\mathcal{X}$.

For each $I \subseteq [n]$, and $x = (x_1,\ldots,x_n) \in \mathcal{X}$ denote $x_I = (x_j)_{j\in I}$, $\bar{x}_I = x_{I^c} = (x_j)_{j\notin I} \in \prod_{j\notin I}\mathcal{X}_i$.
Let $\mu_I$ be the marginal of $\mu$ corresponding to the coordinates indexed by $I$ and $\mu_I(\cdot|\bar{x}_I)$ denote the regular conditional distribution of $x_{I}$ given $\bar{x}_I$ on the probability space $(\mathcal{X},\mu)$. Thus for any Borel set $A\subset \prod_{i\in I} \mathcal{X}_i$, we have
\begin{displaymath}
  \mu(A\times \prod_{j\notin I} \mathcal{X}_j) = \int_{\prod_{j \in I^c} \mathcal{X}_j} \mu_I(A|\bar{x}_I)\mu_{I^c}(d \bar{x}_I) = \int_{\mathcal{X}} \mu_I(A|\bar{x}_I) \mu(dx).
\end{displaymath}
If $I = \{i\}$ we will write e.g., $\bar{x}_i, \mu_i$ instead of $\bar{x}_{\{i\}}$, $\mu_{\{i\}}$.

Recall that for a probability measure $\mu$ and a nonnegative function $f$, the entropy of $f$ relative to $\mu$ is defined as $\Ent_\mu(f) = \E_\mu f\log f - \E_\mu f\log\E_\mu f$ whenever $\E_\mu f\log f < \infty$ and $\Ent_\mu(f) = \infty$ otherwise.

\medskip

The following definition will play a crucial part in what follows.

\begin{defi}[Approximate tensorization of entropy] We will say that $\mu$ has the approximate tensorization property with constant $C$ (abbrev. $AT(C)$) if for every function $f\colon \mathcal{X} \to [0,\infty)$,
\begin{displaymath}
  \Ent_\mu(f) \le C \E_\mu \sum_{i=1}^n\Ent_{\mu_i(\cdot|\bar{x}_i)}(f)
\end{displaymath}
\end{defi}

It is well known that product measures satisfy $AT(1)$, see e.g., \cite[Proposition 5.6]{MR1849347}. Recently Marton \cite{2015arXiv150702803M} (see also \cite{MR3434252,2018arXiv180106348G}) proved the following sufficient condition for tensorization of entropy in discrete product spaces.

\begin{theorem}[Marton, G\"otze-Sambale-Sinulis]\label{thm:Marton-GSS}
  Let $\mu$ be a measure with full support on $\mathcal{X}$. Set
\begin{align}\label{eq:beta}
  \beta = \min_{i\le n}\min_{x\in \mathcal{X}}\mu_i(\{x_i\} |\bar{x}_i)
\end{align}
   Let also $A = (a_{ij})_{i,j\le n}$  satisfy $a_{ii} = 0$ for all $i$ and for $i\neq j$,
\begin{displaymath}
\|\mu_i(\cdot|\bar{x}_i) - \mu_i(\cdot|\bar{y}_i)\|_{TV} \le a_{ij}
\end{displaymath}
for every $x,y\in \mathcal{X}$ which differ only at the $j$-th coordinate.
Assume moreover that $\|A\|_{\ell_2^n\to\ell_2^n} < 1$. Then the measure $\mu$ hast the approximate tensorization property with constant $C = 2\frac{1}{\beta(1-\|A\|_{\ell_2^n\to \ell_2^n})^2}$.
\end{theorem}

In particular in \cite{2018arXiv180106348G} the Authors verify that under our main assumptions the approximate tensorization property is satisfied by the Ising model.

\begin{cor}[G\"otze-Sambale-Sinulis]
If a measure $\mu$ on $\{-1,1\}^n$ is defined by \eqref{eq:Ising} then for any distinct $i,j \in [n]$  and any $x,y \in \mathcal{X}^n$ differing only at the $j$-th coordinate
\begin{displaymath}
\sup_{{x,y \in \mathcal{X}^n}\atop{\bar{x}_j = \bar{y}_j}} \|\mu_i(\cdot|\bar{x}_i) - \mu_i(\cdot|\bar{y}_i)\|_{TV} \le |J_{ij}|
\end{displaymath}
In particular if the condition  \eqref{eq:Dobrushin} is satisfied, then
\begin{displaymath}
\|A\|_{\ell_2^n \to \ell_2^n} \le \|J\|_{\ell_1^n\to \ell_1^n} \le 1 - \rho.
\end{displaymath}
If additionally the vector $h$ satisfies  \eqref{eq:h-bound}, then the coefficient $\beta$ defined in \eqref{eq:beta}, satisfies
\begin{displaymath}
  \frac{1}{C_{\alpha,\rho}} \le \beta \le C_{\alpha,\rho}
\end{displaymath}
where $C_{\alpha,\rho}$ depends only on $\alpha$ and $\rho$.
As a consequence the measure $\mu$ satisfies $AT(C)$ with $C$ depending only on $\rho$ and $\alpha$.
\end{cor}

We will also need the definition of the discrete gradient on $\mathcal{X}$, induced by the measure $\mu$. To this end we will slightly abuse the notation and write $(\bar{x}_i,y_i)$ for the sequence $z$ such that $z_i = y_i$
and $\bar{z}_i = \bar{x}_i$.

Following \cite{2018arXiv180106348G} let us introduce
\begin{defi}
For a measurable function $f \colon \mathcal{X} \to \R$ and $i \in [n]$, define
\begin{align*}
  \mathfrak{d}_i f(x) &= \Big(\frac{1}{2}\int_{\mathcal{X}_i} (f(x) - f(\bar{x}_i,y))^2 \mu_i(dy|\bar{x}_i)\Big)^{1/2}\\
  & = \Big(\frac{1}{2} \E \Big(\Big(f(X) - f(X_1,\ldots,X_{i-1},\widetilde{X}_i,X_{i+1},\ldots,X_n)\Big)^2\Big|X = x\Big)\Big)^{1/2},
\end{align*}
where $X = (X_1,\ldots,X_n)$ is distributed according to $\mu$ and the conditional distribution of $\widetilde{X}_{i}$ given $X=x$ equals to $\mu_i(\cdot|\bar{x}_i)$.

We will denote $\mathfrak{d} f(x) = (\mathfrak{d}_i f(x))_{i\in [n]}$ and regard this vector as an element of $\R^n$ endowed with the standard Euclidean norm $|\cdot|$.
\end{defi}

\begin{defi}
We will say that $\mu$ satisfies the logarithmic Sobolev inequality  with constant $C$ (abbrev. $LSI(C)$) if for every $f\colon \mathcal{X} \to \R$,
\begin{align}\label{eq:log-Sobolev}
  \Ent_\mu (f^2) \le 2C\E |\mathfrak{d} f|^2.
\end{align}
\end{defi}

\begin{remark} The notion of logarithmic Sobolev inequality introduced above can be interpreted as the usual logarithmic Sobolev inequality equivalent to hypercontractivity of the related Glauber dynamics/Gibbs sampler (see e.g., \cite{MR3434252,2015arXiv150702803M,2018arXiv180106348G}), however we will not use this interpretation in the sequel.
\end{remark}

Using the approximate tensorization property together with a log-Sobolev inequality for two point distributions \cite{MR1410112} and a Herbst-type argument of Aida and Stroock \cite{MR1258492} (see also \cite{Bobkov-growth,MR3383337,MR3743923}), G\"otze, Sambale and Sinulis  \cite{2018arXiv180106348G} proved

\begin{theorem}\label{thm:log-Sobolev-Ising}
Let $\mu$ be a measure on $\{-1,1\}^n$, defined by \eqref{eq:Ising}, and assume that $|h|_\infty \le \alpha$ and the Dobrushin condition \eqref{eq:Dobrushin} is satisfied with some $\rho < 1$.  Then there exists a constant $C=C(\alpha,\rho)$ such that $\mu$ satisfies the $LSI(C)$.
As a consequence, if $X$ is a random vector distributed according to $\mu$, then for any $p \ge 2$,
\begin{displaymath}
  \|f(X)\|_p^2 \le \|f(X)\|_2^2 +  2C(p-2)\|\mathfrak{d} f(X)\|_p^2.
\end{displaymath}
\end{theorem}

Combining the above result with the well known fact that the logarithmic Sobolev inequality \eqref{eq:log-Sobolev} implies the Poincar\'e inequality with constant $C$, i.e.,
\begin{displaymath}
  \Var_\mu f \le C \E |\mathfrak{d} f(X)|^2,
\end{displaymath}
we immediately obtain

\begin{cor}\label{cor:Aida-Stroock}
Under the assumptions and notation of Theorem \ref{thm:log-Sobolev-Ising}, for every $f\colon \{-1,1\}^n \to \R$, and any $p \ge 2$,
\begin{displaymath}
  \|f(X) - \E f(X)\|_p \le \sqrt{2Cp}\|\mathfrak{d} f(X)\|_p.
\end{displaymath}

\end{cor}

This inequality will be the basis of our inductive argument in the proof of Theorem \ref{thm:polynomial-Ising}.

\section{Proof of Theorem \ref{thm:polynomial-Ising}}\label{sec:proof-Ising}
\begin{proof}[Proof of Theorem \ref{thm:polynomial-Ising}]
The proof will be an adaptation of arguments from \cite{MR3383337} (see also \cite{MR3743923}) to the discrete case.
To carry it out it will be convenient to introduce an inner product $\langle \cdot,\cdot \rangle$ on the space of $k$-tensors with the formula
\begin{displaymath}
  \langle A,B\rangle = \sum_{|\ii|\le n} a_\ii b_\ii,
\end{displaymath}
where $A = (a_\ii)_{|\ii|\le n}$, $B = (b_\ii)_{|\ii|\le n}$. Let us also recall the notation $x^1\otimes \cdots\otimes x^k = (x^1_{i_1}x^2_{i_2}\cdots x^{k}_{i_k})_{i_1,\ldots,i_k\le n}$ for any vectors $x^j = (x^j_1,\ldots,x^j_n)$, $j=1,\ldots,k$.

We will first prove by induction on $d$ that for any positive integer $d$, and any function $f\colon \{-1,1\}^n \to \R$,
\begin{align}\label{eq:moment-estimate}
  \|f(X) - \E f(X)\|_p \le K^{d}\|\langle \nabla^d f(X),G_1\otimes\cdots\otimes  G_d\rangle\|_p + \sum_{i=1}^{d-1}  K^{i}\|\langle \E_X\nabla^i f(X),G_1\otimes\cdots\otimes  G_i\rangle\|_p,
\end{align}
where $K$ is a constant depending only on $\alpha,\rho$ and $G_1,\ldots,G_d$ are i.i.d.\ standard Gaussian vectors in $\R^n$ independent of $X$. Here by $\E_X$ we denote expectation with respect to the random vector $X$ and (as explained in Section \ref{sec:main-results}) the derivatives $\nabla^i f$ denote the derivatives of the tetrahedral polynomial coming from the Fourier--Walsh expansion of $f$. We remark that $d$ does not necessarily coincide with the degree of $f$.

To this end we will proceed by induction. Consider thus $f(x) = A_0 + \sum_{k=1}^D \langle A^k,x^{\otimes k}\rangle$, where $A_0\in \R$ and for $k=1,\ldots,D$, $A^k = (a^k_{\ii})_{|\ii|\le n}$ are $k$-indexed symmetric matrices with zeros on generalized diagonals.

Let $\widetilde{X}_i$, $i=1,\ldots,n$ be $\{-1,1\}$-valued random variables (possibly defined on some extension of the original probability space) such that the conditional distribution of
$\widetilde{X}_i$ given $X=x$ equals $\mu_i(\cdot|\bar{x}_i)$.

Recall that for $i\le n$,
\begin{align*}
 \mathfrak{d}_i f(x) =  \Big(\frac{1}{2} \E \Big(\Big(f(X) - f(X_1,\ldots,X_{i-1},\widetilde{X}_i,X_{i+1},\ldots,X_n)\Big)^2\Big|X=x\Big)\Big)^{1/2},
\end{align*}
hence by Corollary \ref{cor:Aida-Stroock}, using the notation $\bar{X}_i = (X_j)_{j\neq i}$ we have for any $p\ge 2$,
\begin{displaymath}
  \|f(X)-\E f(X)\|_p \le \sqrt{Cp}\Big(\E \Big|\sum_{i=1}^n \E \Big((f(X) - f(\bar{X}_i,\widetilde{X}_i))^2|X\Big)\Big|^{p/2}\Big)^{1/p}.
\end{displaymath}
Using Jensen's inequality for the conditional expectation, we can further write
\begin{displaymath}
\|f(X)-\E f(X)\|_p \le \sqrt{Cp}\Big(\E \Big|\sum_{i=1}^n (f(X) - f(\bar{X}_i,\widetilde{X}_i))^2\Big|^{p/2}\Big)^{1/p}.
\end{displaymath}

Define now for $k=1,\ldots,d$ and $i=1,\ldots,n$ the $(k-1)$-indexed matrices $A^{k,i} = (a^{k,i}_{\ii})_{|\ii|\le n} = (a^k_{(i,\ii)})_{|\ii|\le n}$, where for $\ii = (i_1,\ldots,i_{k-1})$ we write $(i,\ii) = (i,i_1,\ldots,i_{k-1})$.
Using the fact that the generalized diagonals of the matrices $A^k$ vanish together with the symmetry of $A^k$, we get
\begin{align*}
f(X) - f(\bar{X}_i,\widetilde{X}_i) &= (X_i - \widetilde{X}_i) \sum_{k=1}^D k\langle A^{k,i}, X^{\otimes (k-1)}\rangle = (X_i-\widetilde{X}_i)\frac{\partial}{\partial x_i} f(X).
\end{align*}
Since $|X_i - \widetilde{X}_i| \le 2$, by combining this equality with the previous estimate, we obtain
\begin{displaymath}
  \|f(X)-\E f(X)\|_p \le 2 \sqrt{Cp} (\E |\nabla f(X)|^p)^{1/p}.
\end{displaymath}

Using the fact that if $g$ is a standard Gaussian variable, then for $p \ge 1$ we have $\sqrt{p}M^{-1}\le  \|g\|_p \le M \sqrt{p}$, where $M$ is a universal constant,
we can write the above inequality as
\begin{align}\label{eq:moment-estimate-d=1}
  \|f(X)-\E f(X)\|_p \le K \Big\|\Big\langle \nabla f(X),G\Big\rangle \Big\|_p
\end{align}
for a standard $n$-dimensional Gaussian vector $G$, independent of $X,\{\widetilde{X}_i\}_{i\le n}$ and $K = 2\sqrt{C}M$.  This establishes \eqref{eq:moment-estimate} for $d=1$.

The induction step follows just by the case $d=1$ and the triangle inequality in $L_p$. Indeed, assuming that \eqref{eq:moment-estimate} holds for $d$, by the triangle inequality and linearity of expectation we get
\begin{align*}
  \|f(X) - \E f(X)\|_p \le &K^{d}\|\langle \nabla^d f(X),G_1\otimes\cdots\otimes  G_d\rangle - \E_X \langle \nabla^d f(X),G_1\otimes\cdots\otimes  G_d\rangle\|_p \\
  &+ \sum_{i=1}^{d}  K^{i}\|\langle \E_X\nabla^i f(X),G_1\otimes\cdots\otimes  G_i\rangle\|_p,
\end{align*}
Applying now (conditionally on $G_1,\ldots,G_d$) \eqref{eq:moment-estimate-d=1}  to the first term on the right hand side and using the Fubini theorem we obtain
\begin{displaymath}
  \|\langle \nabla^d f(X),G_1\otimes\cdots\otimes  G_d\rangle - \E_X \langle \nabla^d f(X),G_1\otimes\cdots\otimes  G_d\rangle\|_p \le K \|\langle \nabla^{d+1} f(X),G_1\otimes\cdots\otimes G_{d+1}\rangle\|_p
\end{displaymath}
(note that $\langle \nabla^d f(X),G_1\otimes\cdots \otimes G_d\rangle$ is tetrahedral as a polynomial in $X$).

This ends the induction step and establishes \eqref{eq:moment-estimate}.

\medskip

If $f$ is a polynomial of degree $d$, then $\nabla^d f(X)$ is deterministic (and thus equal to its expectation) so \eqref{eq:moment-estimate} can be written in a more concise way
\begin{displaymath}
   \|f(X)-\E f(X)\|_p \le \sum_{i=1}^{d}  K^{i}\|\langle \E_X\nabla^i f(X),G_1\otimes\cdots\otimes  G_i\rangle\|_p
\end{displaymath}

We will now use a result by Lata{\l}a \cite{MR2294983}, which asserts the existence of constants $C_k$, depending only on $k$, such that for any $k$-index matrix $A$, and $p\ge 2$,
\begin{displaymath}
  \|\langle A,G_1\otimes\cdots\otimes G_k\rangle\|_p\le C_k\sum_{\mathcal{I}\in P_k} p^{|\mathcal{I}|/2}\|A\|_\mathcal{I},
\end{displaymath}
which yields
\begin{displaymath}
  \|f(X)-\E f(X)\|_p \le  C \sum_{k=1}^d\sum_{\mathcal{I}\in P_k}  p^{|\mathcal{I}|/2} \|\E\nabla^k f(X)\|_\mathcal{I},
\end{displaymath}
where $C$ depends on $\rho,\alpha,d$.

By Chebyshev's inequality in $L_p$ this gives for $p\ge0$,
\begin{displaymath}
  \p\Big(|f(X)-\E f(X)|\ge Ce \sum_{k=1}^d\sum_{\mathcal{I}\in P_k}  p^{|\mathcal{I}|/2} \|E\nabla^k f(X)\|_\mathcal{I} \Big) \le e^{2-p}
\end{displaymath}
(the additional factor $e^2$ on the right hand side allows to extend the estimate from $p\ge 2$ to all $p \ge 0$).
The theorem follows now by a change of variables and adjustment of constants.
\end{proof}

\section{Convex concentration and improved estimates for positive definite quadratic forms}\label{sec:improved-concentration}

The estimates of Theorem \ref{thm:polynomial-Ising} are of Gaussian nature, i.e., they have the same form as two-sided estimates valid for polynomials in independent Gaussian variables. Since the values of the random variables in the Ising model are $\pm 1$, it is natural to look for estimates resembling those known for polynomials in independent Rademacher variables. In this case the situation is however more complicated, as two-sided bounds are known only for polynomials of degree at most 3 (see \cite{MR1338834,MR1686370,MR3052405}).

Below in Theorem \ref{thm:quadratic-forms-improved} we present estimates similar in nature to those for Rademacher sequences for quadratic forms $\langle A X,X\rangle$ where $A$ is a non-negative definite matrix and $X$ is a random vector with bounded coefficients, satisfying the approximate tensorization property. In some situations they improve on the bounds one can get for the Ising model from Theorem \ref{thm:polynomial-Ising}, however in general they are not comparable to them, because they involve norms of the matrix $A$ and not just its off-diagonal part (note that in the case of the Ising model, the contribution from the diagonal is deterministic). It is natural to conjecture that (similarly as for the i.i.d.\ case) the assumption of non-negative definiteness is an artefact of our proof and can be actually dropped, however at present we are not able to obtain such more general bounds.

\subsection{Convex concentration}

As a tool for proving estimates for quadratic forms we will derive concentration inequalities for convex Lipschitz functions for measures on products of compact sets, satisfying the approximate tensorization property, which are of independent interest. In particular this will allow us to obtain concentration for linear combinations with vector coefficients  (see Proposition \ref{prop:vector-valued}), which generalize moment estimates obtained in the Rademacher case by Dilworth and Montgomery-Smith \cite{dilworth1993} (see also \cite{MR1686370}).

Recall that a random vector $X$ in $\R^n$ has the convex concentration property with constant $K$ if for any $L$-Lipschitz convex function $f\colon \R^n \to \R$, and any $t > 0$,
\begin{align}\label{eq:convex-concentration-median}
  \p(|f(X) - \Med f(X)| \ge t) \le 2\exp(-t^2/K^2 L^2).
\end{align}
It is well known that the above property is up to constant equivalent to concentration around the mean, i.e.,
\begin{align}\label{eq:convex-concentration-mean}
  \p(|f(X) - \E f(X)| \ge t) \le 2\exp(-t^2/\widetilde{K}^2 L^2),
\end{align}
i.e., the inequalities \eqref{eq:convex-concentration-median} and \eqref{eq:convex-concentration-mean} imply each other and the constants $K$ and $\widetilde{K}$ depend only on one another.
We will say that a random vector $Z$ satisfies the dimension-free convex concentration property with constant $K$ if for any $N$, the random vector $X = (X_1,\ldots,X_N)$, where $X_i$ are i.i.d.\ copies of $Z$, satisfies \eqref{eq:convex-concentration-median}.

We will now relate the approximate tensorization property of measures on $[-1,1]^n$ to the convex concentration property, showing in particular that if the distribution of $X$ is given by \eqref{eq:Ising}, where $|h|_\infty \le \alpha$ and $J_{ij}$ satisfy the Dobrushin condition \eqref{eq:Dobrushin}, then $X$ satisfies the dimension-free convex concentration property with a constant depending only on $\alpha$ and $\rho$ (Proposition \ref{prop:from-AT-to-convex-concentration} below). Next we will prove that convex concentration property for measures on products of compact sets can be in fact improved by taking into account the uniform bounds on the components of $X$. Finally we will illustrate this phenomenon with applications to linear forms with vector coefficients and quadratic non-negatively definite forms (Proposition \ref{prop:vector-valued} and Theorem \ref{thm:quadratic-forms-improved}).

In order to pass from approximate tensorization of entropy to dimension-free convex concentration property, we will use weak transportation inequalities, introduced recently by Gozlan, Roberto, Samson and Tetali \cite{MR3706606}.

Let us denote by $\mathcal{P}_1(\R^n)$ the set of all probability measures on $\R^n$ with finite first moment.

\begin{defi}
  Let $\mu$ and $\nu$ be probability measures on $\R^n$. Assume that $\nu \in \mathcal{P}_1(\R^n)$. For a convex, lower semicontinuous function $\theta\colon \R^n \to [0,\infty]$, such that $\theta(0) = 0$ define the weak transport cost between $\mu$ and $\nu$ as
  \begin{displaymath}
    \cTbar_{\theta}(\nu|\mu) = \inf_{\pi}\int_{\R^n}\theta\Big(x - \int_{\R^n} y p_x(dy)\Big)\mu(dx),
  \end{displaymath}
  where the infimum is taken over all couplings $\pi$ between $\mu$ and $\nu$ (i.e., measures on $\R^n\times \R^n$ with marginals $\mu,
  \nu$) and for $x \in \R^n$, $p_x(\cdot)$ is the conditional measure defined ($\mu$ almost surely) by $\pi(dxdy) = p_x(dy)\mu(dx)$.
\end{defi}

Using probabilistic notation one can write
\begin{align}\label{eq:Tbar-probabilistically}
  \cTbar_{\theta}(\nu|\mu) = \inf_{(X,Y)} \E \theta(X - \E(Y|X)),
\end{align}
where the infimum is taken over all pairs of random vectors $(X,Y)$ with values in $\R^n \times \R^n$, such that $X$ is distributed according to $\mu$ and $Y$ according to $\nu$.

Recall also that if $\mu,\nu$ are two probability measures then the relative entropy of $\nu$ with respect to $\mu$ is given by the formula
\begin{displaymath}
  H(\nu|\mu) = \E_\nu \log\Big(\frac{d\nu}{d\mu}\Big)
\end{displaymath}
if $\nu$ is absolutely continuous with respect to $\mu$ and $H(\nu|\mu) = \infty$ otherwise.

\begin{defi} Let $\mu \in \mathcal{P}_1(\R^n)$ and $\theta\colon \R^n \to [0,\infty]$ be a~convex lower semicontinuous function with $\theta(0) = 0$. We will say that $\mu$ satisfies the inequality
 $\Tbar_{\theta}$ if for every probability measure $\nu\in \mathcal{P}_1(\R^n)$,
 \begin{align}\label{eq-Tbar}
	\max\Big(\cTbar_{\theta}(\nu|\mu),\cTbar_{\theta}(\mu|\nu)\Big) \le H(\nu|\mu).
\end{align}
\end{defi}

The following theorem established in \cite{MR3706606} describes connections between dimension-free convex  concentration, weak transportation inequalities and log-Sobolev inequalities for convex and concave functions.
\begin{theorem}\label{thm:GRST}
Let $X$ be a random vector in $\R^n$ with distribution $\mu$.
The following conditions are equivalent.
\begin{enumerate}[(i)]
  \item There exists $K$ such that $X$ has the dimension-free convex concentration property with constant $K$.
  \item There exists $c$ such that $\mu$ satisfies the inequality $\Tbar_\theta$ with $\theta(x) = c|x|^2$.
  \item There exist $D,\lambda > 0$ such that for every convex Lipschitz function and every concave function whose Hessian is bounded from below by $(-\lambda){\rm Id}$,
  \begin{align}\label{eq:convex-concave-LSI}
    \Ent e^{f(X)} \le D \E |\nabla f(X)|^2 e^{f(X)}.
  \end{align}
\end{enumerate}
Moreover for any two assertions above the constants in one of them may be taken to depend only on the constants in the other one.
\end{theorem}

Using the above result we can easily obtain the following proposition, which may be useful e.g., in statistical applications, when dealing with i.i.d.\ samples drawn from the measure $\mu$ (see e.g., \cite{MR3775918} for a discussion of applications related to the Ising model).

\begin{prop}\label{prop:from-AT-to-convex-concentration}
If $X$ is a $[-1,1]^n$-valued random vector with law $\mu$, which satisfies the approximate tensorization $AT(C)$, then $X$ satisfies the  dimension-free convex concentration inequality with constant $K$, depending only on $C$.
\end{prop}

\begin{proof}
The celebrated convex distance inequality by Talagrand (see eg. \cite{MR1361756,MR1419006}) asserts that any random variable with support in $[-1,1]$ satisfies the dimension-free convex concentration property with a universal constant. In particular by Theorem \ref{thm:GRST} it satisfies the log-Sobolev inequality \eqref{eq:convex-concave-LSI} with some universal constants $D$, $\lambda$. Consider any function $f\colon \R^n \to \R$, which is either convex or concave with $\nabla^2 f(x) \ge -\lambda {\rm Id}$ for all $x$. In particular in the latter case for any $i \le n$, $\frac{\partial^2 f(x)}{\partial x_i^2} \ge -\lambda$. One can thus apply the one-dimensional version of \eqref{eq:convex-concave-LSI} to $\mu_i(\cdot|\bar{x}_i)$ and the function $x_i \mapsto f(\bar{x}_i,x_i)$, which together with the condition $AT(C)$ gives \eqref{eq:convex-concave-LSI} with constants $CD$ and $\lambda$. The proof is now concluded by another application of Theorem \ref{thm:GRST}.
\end{proof}

It is easy to see that if a measure $\mu$ supported on $[-1,1]^n$ satisfies $\Tbar_\theta$ with $\theta(x) = c|x|^2$ then it actually satisfies a stronger inequality
$\Tbar_{\gamma}$ with $\gamma(x) = |x|^2$ if $|x|_\infty < 2$ and $\gamma(x) = \infty$ otherwise. Indeed for the right-hand side to be finite $\nu$ must be also supported on $[-1,1]^n$, in which case by \eqref{eq:Tbar-probabilistically} $\cTbar_\theta(\mu|\nu) = \cTbar_\gamma(\mu|\nu)$ and $\cTbar_\theta(\nu|\mu) = \cTbar_\gamma(\nu|\mu)$.
In fact, weak transportation inequalities with such strengthened cost functions can hold only for compactly supported measures (see  \cite{2017arXiv170204698S}). The interest in such strengthening lies in the fact that by taking into account the boundedness of random variables, it implies concentration inequalities  stronger than the subgaussian  bound given by \eqref{eq:convex-concentration-mean} (see e.g., \cite{2017arXiv170301765A} for concentration results corresponding to various cost functions $\theta$). As shown in the next proposition, such inequalities can be also easily inferred just at the level of convex concentration. To formulate this result let us introduce a family of norms on $\R^n$ given for $p > 0$ by the formula
\begin{align}\label{eq:Montgomery-Smith-norm}
  \|x\|_{\{1\},p} = \sup\Big\{\sum_{i=1}^n x_iy_i\colon |y| \le \sqrt{p}, |y|_\infty \le 1\Big\}.
\end{align}

It is not difficult to see that
\begin{align}\label{eq:dual-representation}
  \|x\|_{\{1\},p} & \le  \sum_{i\le p} x_i^\downarrow + \sqrt{p}\sqrt{\sum_{i >p} (x_i^\downarrow)^2} \le 2 \|x\|_{\{1\},p},
\end{align}
where $(x^\downarrow_i)_{i\le n}$ is the nonincreasing rearrangement of the sequence $(|x_i|)_{i\le n}$. In fact one has
\begin{align}\label{eq:partial-dual}
  \max\Big(\sum_{i\le p} x_i^\downarrow,  \sqrt{p}\sqrt{\sum_{i >p} (x_i^\downarrow)^2}\Big) \le \|x\|_{\{1\},p}.
\end{align}

Such norms are equivalent to interpolation norms between the spaces $\ell_2^n$ and $\ell_1^n$ and in a probabilistic context appeared for the first time in the paper \cite{MR1013975}, where it is shown that if $\varepsilon_1,\ldots,\varepsilon_n$ are independent Rademacher variables, then for $x\in \R^n$ and $p\ge 2$,
\begin{displaymath}
  \frac{1}{C} \|x\|_{\{1\},p} \le \Big\|\sum_{i=1}^n \varepsilon_i x_i\Big\|_p \le C\|x\|_{\{1\},p}
\end{displaymath}
where $C$ is a universal constant.
The meaning of the subscript $\{1\}$ will become clear when we define counterparts of this norm for matrices. To keep uniform notation, we introduce it already here.

\medskip

We are now ready to state the strengthened concentration result.

\begin{prop}\label{prop:vector-valued}
Let $X$ be a random vector with values in $[-1,1]^n$, satisfying \eqref{eq:convex-concentration-median}. Then for any smooth convex Lipschitz function $f\colon \R^n \to \R$
and any $p > 0$,
\begin{align}\label{eq:improved-convex-concentration}
  \p(|f(X) - M| > C \sup_{x} |\nabla f(x)|_{\{1\},p}) \le 4\exp(-p/K^2).
\end{align}
where $M$ is the mean or the median of $f(X)$ and $C$ is a universal constant.
\end{prop}

This improves on what can be obtained from \eqref{eq:convex-concentration-mean} since as one can see from \eqref{eq:dual-representation}, $\|x\|_{\{1\},p} \le \sqrt{p}|x|$.

\begin{remark}\label{rem:general-function}  By standard regularization arguments (see e.g., see \cite[p.~429]{MR1756011}) one can drop the smoothness assumptions on $f$, by replacing $\sup_x \|\nabla f(x)\|_{\{1\},p}$ with the Lipschitz constant of $f$ with respect to the norm dual to $\|\cdot\|_{\{1\},p}$. One can also assume that $f$ is defined on $[-1,1]^n$ since one can extend it to $\R^n$ with the formula $\tilde{f}(y) = \sup_{x\in(-1,1)^n} (f(x) + \langle \nabla f(x),y-x\rangle)$, without altering the Lipschitz constant (here $\nabla f(x)$ denotes some subgradient of $f$ at $x$).
\end{remark}

Before proving the above proposition, let us illustrate it with examples, to show how it improves on the usual subgaussian convex concentration \eqref{eq:convex-concentration-mean}.

In view of Remark \ref{rem:general-function}, Proposition \ref{prop:vector-valued} yields the following corollary.

\begin{cor}\label{cor:vector}
Let $X = (X_1,\ldots,X_n)$ be a $[-1,1]^n$-valued random vector satisfying the convex concentration property \eqref{eq:convex-concentration-median}. Let $(E,\|\cdot\|)$ be a Banach space with dual $(E^\ast,\|\cdot\|_\ast)$, $v_1,\ldots,v_b \in E$ and let
\begin{displaymath}
  Z = \Big\|\sum_{i=1}^n v_i X_i\Big\|,
\end{displaymath}
then for every $p > 0$,
\begin{align}\label{eq:vector}
  \p(|Z - \E Z| \ge \sup_{\varphi \in E^\ast, \|\varphi\|_\ast \le 1} \|(\varphi(v_i))_{i\le n}\|_{\{1\},p}) \le 4 e^{-p/K^2}.
\end{align}
\end{cor}

\begin{remark}
In \cite{dilworth1993}, Dilworth and Montgomery-Smith proved that if $X$ is a vector of independent Rademacher variables, then
$\p(Z \ge 2\E Z + \sup_{\varphi \in E^\ast, \|\varphi\|_\ast \le 1} \|(\varphi(v_i))_{i\le n}\|_{\{1\},p}) \le 2 e^{-cp}$.
The fact that one can improve this statement to concentration around mean, seems to have become by then a part of the folklore.
\end{remark}

\begin{example}\label{ex:easy-sum}
Let us now provide a simple one dimensional example. For illustration purposes it will be more convenient to state it in terms of infinite sequences of random variables. Let thus $X_1,X_2,\ldots$ be centered random variables with values in $[-1,1]$, such that for all $n$, the vector $X = (X_1,\ldots,X_n)$ satisfies \eqref{eq:convex-concentration-median} with $K=1$ (for simplicity). Define the random variable $Z = \sum_{i=1}^\infty \frac{1}{i}X_i$ (note that thanks  to \eqref{eq:convex-concentration-mean} this sequence converges in $L_2$). Then by \eqref{eq:convex-concentration-mean}, we get $\p(|Z| \ge t) \le 2\exp(-ct^2)$ for some $c>0$. However, it is easy to see that $\|(1/i)_{i=1}^\infty\|_{\{1\},p} \simeq \log p$ (up to multiplicative constants) for $p \ge 2$ thus by Proposition \ref{prop:vector-valued} (after adjusting the constants) we obtain
\begin{displaymath}
  \p(|Z| \ge t) \le 4e^{-e^{c't}}.
\end{displaymath}
This Gumbel type tail decay is clearly much faster than Gaussian. Note that thanks to the one-dimensional nature this example can be in fact easily recovered directly from \eqref{eq:convex-concentration-mean} by combining it with obvious pointwise bounds on the variables $X_i$. Since \eqref{eq:vector} can be equivalently restated as
\begin{displaymath}
  \p\Big(|Z - \E Z| \ge \sup_{x \in T} \|x\|_{\{1\},p}\Big) \le 4e^{-p/K^2}
\end{displaymath}
for any bounded set $T$ of vectors and $Z = \sup_{x\in T} |\sum_i x_iX_i|$, one can easily create more complicated examples with various types of tail decay.
\end{example}



\medskip

\begin{proof}[Proof of Proposition \ref{prop:vector-valued}]
The idea of the proof goes back to Talagrand and is by now classical. The main additional observation one needs to make is that exploring  the boundedness of the support may lead to improved inequalities for general convex Lipschitz functions rather than just for linear functions.

We will start by proving the inequality in question with the median. Let thus $M = \Med f$ and consider first the convex set $A = \{x\in[-1,1]^n\colon f(x) \le M\}$, so that $\p(X \in A) \ge 1/2$.
Define $g(x) = {\rm dist}(x,A)$, then $\Med g(X) = 0$ and by convexity of $A$, $g$ is a convex function. Note that if for $z \in [-1,1]^n$, $f(z) >  M + 3 \sup_{x} |\nabla f(x)|_{\{1\},p}$, then by convexity for any $y \in A$,
\begin{align*}
  M + 3\|\nabla f(z)\|_{\{1\},p} & < f(z) \le f(y) + \langle z - y,\nabla f(z)\rangle \\
  &\le M + 2 \sum_{i\le p} (\frac{\partial}{\partial x_i} f(z))^\downarrow + \frac{1}{\sqrt{p}}|z-y| \sqrt{p}\sqrt{\sum_{i>p}((\frac{\partial}{\partial x_i} f(z))^\downarrow)^2}\\
  & \le M + 2\|\nabla f(z)\|_{\{1\},p} + \frac{1}{\sqrt{p}} |z-y| \|\nabla f(z)\|_{\{1\},p},
\end{align*}
(where we used \eqref{eq:partial-dual}) and so $|z - y| > \sqrt{p}$. Taking infimum over all $y \in A$, and recalling that $\Med g(X) = 0$, we obtain
\begin{align}\label{eq:bla-bla-upper-tail}
  \p(f(X) > M + 3 \sup_{x} |\nabla f(x)|_{\{1\},p}) \le \p(g(X) \ge \Med g(X) + \sqrt{p}) \le 2\exp(-p/K^2),
\end{align}
where in the last inequality we used \eqref{eq:convex-concentration-median}.

As for the lower tail, we can clearly assume that $f$ is not constant. In particular $\sup_x \|\nabla f(x)\|_{\{1\},p} > 0$. Then, denoting
\begin{displaymath}
A = \{x\in[-1,1]^n\colon f(x) \le M - 3 \sup_x \|\nabla f(x)\|_{\{1\},p}\},
\end{displaymath}
by similar estimates as above one obtains that for $z \in [-1,1]^n$, if ${\rm dist}(z,A) < \sqrt{p}$, then $f(z) < M$. Thus, denoting $B = \{x\in \R^n\colon f(x) \ge M\}$ we get $A \subseteq \{x\in [-1,1]^n\colon {\rm dist}(x,B) \ge \sqrt{p}\}$. Since $A \subseteq B^c$ we can also assume that $B^c\neq \emptyset$. The function $g(x) = {\rm dist}(x,B)$ is 1-Lipschitz, concave on the complement of $B$ and can be extended to a function $\widetilde{g}(x) := \inf_{z\in B^c} (g(z) + \langle \nabla g(z),x-z\rangle)$, which is 1-Lipschitz, concave on $\R^n$ and non-positive on $B$. Moreover $\widetilde{g} = g$ on $B^c$.
Thus $\Med \widetilde{g}(X) \le 0$ and so
\begin{displaymath}
A \subseteq \{x\in [-1,1]^n\colon g(x) \ge \sqrt{p}\} \subseteq \{x \in [-1,1]^n \colon \widetilde{g}(x) \ge \Med \widetilde{g}(X) + \sqrt{p}\}.
\end{displaymath}
 As a consequence $\p(X \in A) \le 2\exp(-p/K^2)$, which together with \eqref{eq:bla-bla-upper-tail} proves that
 \begin{displaymath}
   \p\Big(|f(X) - \Med f(X)| \ge 3\sup_x \|\nabla f(x)\|_{\{1\},p}) \le 4e^{-p/K^2}.
 \end{displaymath}

To pass from the median to the mean, we notice that for $t \ge 1$, $\|x\|_{\{1\},t p} \le \sqrt{t}\|x\|_{\{1\},p}$, so applying the above estimate with $t^2p$ instead of $p$, we get
\begin{displaymath}
  \p(|f(X)- \Med f(X)| \ge 3t \sup_{x} \|\nabla f(x)\|_{\{1\},p}) \le 4\exp(-t^2p/K^2).
\end{displaymath}
In particular by Jensen's inequality and integration by parts this yields
\begin{displaymath}
|\E f(X) - \Med f(X)| \le \E|f(X) - \Med f(X)| \le (3+C'Kp^{-1/2})\sup_{x} \|\nabla f(x)\|_{\{1\},p}
\end{displaymath}
for some universal constant $C'$, which gives
\begin{displaymath}
  \p(|f(X) - \E f(X)| \ge (6 + C'Kp^{-1/2})\sup_{x} \|\nabla f(x)\|_{\{1\},p}) \le 4\exp(-p/K^2).
\end{displaymath}

If $p > K^2$, this gives \eqref{eq:improved-convex-concentration} for $M = \E f(X)$ with  $C = 6 + C'$, otherwise \eqref{eq:improved-convex-concentration} is trivial, as the right hand side exceeds one.
\end{proof}

\subsection{Quadratic forms}
We will now pass to quadratic forms. In order to formulate tail estimates in this case, we need to introduce two additional norms of a symmetric matrix. Following \cite{MR1686370} we define
  \begin{align*}
    \|A\|_{\{1,2\},p} & = \|(A_i)_{i\le n}\|_{\{1\},p},
  \end{align*}
where $A_i = \sqrt{\sum_{j=1}^n a_{ij}^2}$ and
\begin{align*}
\|A\|_{\{1\}\{2\},p} = \sup\Big\{\sum_{i,j=1}^n a_{ij}x_iy_j\colon |x|,|y| \le \sqrt{p},  \;|x|_\infty, |y|_\infty \le 1\Big\}.
\end{align*}

We note that by \eqref{eq:dual-representation},
\begin{displaymath}
   \|A\|_{\{1,2\},p} \le \sum_{i\le p} A_i^\downarrow + \sqrt{p}\Big(\sum_{i>p} (A_i^\downarrow)^2)^{1/2} \le 2 \|A\|_{\{1,2\},p},
\end{displaymath}
which gives a simpler expression if one is interested in concentration up to dimension-free constants.

\begin{remark} It is easy to see that
\begin{displaymath}
  \|A\|_{\{1,2\},p} = \sup\Big\{\sum_{i,j=1}^n a_{ij}x_{ij} \colon \quad \|x\|_{\{1,2\}} \le \sqrt{p},\; \max_{i\le n} \sum_{j=1}^n x_{ij}^2 \le 1\Big\},
\end{displaymath}
which justifies the subscript $\{1,2\}$ used in the notation. It is also clear that $\|A\|_{\{1,2\},p} \le C\sqrt{p}\|A\|_{\{1,2\}} = C\sqrt{p}\|A\|_{HS}$ and $\|A\|_{\{1\}\{2\},p} \le p \|A\|_{\{1\}\{2\}} = p \|A\|_{\ell_2^n\to\ell_2^n}$.
\end{remark}

\medskip

In \cite{MR1686370} Lata{\l}a proved that there exists a universal constant $C$ such that if $X=(\varepsilon_1,\ldots,\varepsilon_n)$, where $\varepsilon_i$'s are independent Rademacher variables, then for any symmetric matrix with vanishing diagonal and any $p\ge 2$, one has
\begin{displaymath}
  C^{-1}(\|A\|_{\{1,2\},p} + \|A\|_{\{1\}\{2\},p}) \le \|\langle AX, X\rangle\|_p \le C(\|A\|_{\{1,2\},p} + \|A\|_{\{1\}\{2\},p}),
\end{displaymath}
similar bounds were also obtained for cubic forms in \cite{MR3052405}. As a consequence,  by Chebyshev's and Paley--Zygmund inequalities (see \cite{MR1686370}), in this case for any $p > 0$,
\begin{displaymath}
\p\Big(|\langle AX, X\rangle| \ge Ce(\|A\|_{\{1,2\},p} + \|A\|_{\{1\}\{2\},p})\Big) \le e^{2-p}
\end{displaymath}
and
\begin{displaymath}
\p\Big(|\langle AX, X\rangle| \ge c(\|A\|_{\{1,2\},p} + \|A\|_{\{1\}\{2\},p})\Big) \ge \min(c,e^{-p}).
\end{displaymath}

The result we obtain for quadratic forms in dependent random variables is

\begin{theorem}\label{thm:quadratic-forms-improved}
  Let $X$ be a centered random vector with values in $[-1,1]^n$, satisfying the convex concentration property with constant $K$ and let $A = (a_{ij})_{i,j\le n}$ be a symmetric nonnegative definite matrix. Then there exists a constant $C_K$, depending only on $K$, such that for any $p \ge 0$,
  \begin{align}\label{eq:upper-estimates-tails}
  \p\Big(\langle AX,X\rangle - \E\langle AX,X\rangle \ge C_K(\|A\|_{\{1,2\},p} + \|A\|_{\{1\}\{2\},p})\Big) \le 4e^{-p}
  \end{align}
  and
\begin{align}\label{eq:lower-estimates-tails}
\p\Big(\langle AX,X\rangle - \E\langle AX,X\rangle \le - C_K\min(\|A\|_{\{1,2\},p} + \|A\|_{\{1\}\{2\},p},\sqrt{p}\|A\|_{\{1,2\}})  \Big) \le 4e^{-p}.
\end{align}
\end{theorem}

\begin{remark}
The mean zero assumption in the above theorem is introduced only to simplify its formulation. Clearly in the general case one can recenter the vector and handle the linear correction by Proposition \ref{prop:vector-valued}.
\end{remark}

\begin{remark}\label{rem:decomposition-fails}
For the Ising model the norms $\|A\|_{\{1,2\},p}$ and $\|A\|_{\{1\}\{2\},p}$ introduce unnecessary contribution from the diagonal of $A$, which does not influence the value of $\langle AX,X\rangle - \E \langle AX,X\rangle$. However for random variables not supported on $\{-1,1\}$, in general this contribution has to be taken into account. It is not difficult to see that up to constants it corresponds to the $\|\cdot\|_{\{1\},p}$ norm of the vector consisting of diagonal elements from $A$, which is consistent with estimates of Proposition \ref{prop:vector-valued} as well as tail bounds for sums of independent bounded random variables.

As already mentioned at the beginning of the section, one expects that the assumption of nonnegative definiteness of the matrix $A$ is not needed in Theorem \ref{thm:quadratic-forms-improved}. In \cite{MR3407216} it is shown that the convex concentration property \eqref{eq:convex-concentration-mean} implies the Hanson-Wright inequality \eqref{eq:Hanson-Wright-123} (with $c$ depending on $K$) for arbitrary matrices by splitting the matrix into the sum of its positive and negative definite parts and treating each of them separately (using convexity). This strategy does not work here, since the $\|\cdot\|_{\mathcal{I},p}$ norms are not invariant under conjugation and the norms of positive and negative parts can be of greater order than the corresponding norms of the original matrix. This can be seen e.g., with a matrix $A = (a_{ij})_{i,j\le n}$ such that $a_{1i} = a_{i1} = 1$ for $i \neq 1$ and all the other coefficients are zero. In this case for $1\ll p \ll n$ we get $\|A\|_{\{1,2\},p} + \|A\|_{\{1\}\{2\},p} \simeq \sqrt{p}\sqrt{n}$, whereas if $A_\pm$ is the positive/negative part of $A$, then $\|A_\pm\|_{\{1,2\},p} + \|A_\pm|_{\{1\}\{2\},p} \simeq p\sqrt{n}$.

\medskip

For the Ising model one may hope that the assumption of nonnegative definiteness of the matrix $A$ in Theorem \ref{thm:quadratic-forms-improved} could be removed by a repetition of the proof of Theorem \ref{thm:polynomial-Ising} with auxiliary Rademacher variables, instead of Gaussian ones, i.e., by proving that for every $f$ (seen as a tetrahedral polynomial) and
$p\ge 2$, $$\|f(X) - \E f(X)\|_p \le C\Big(\E \|\nabla f(X)\|_{\{1\},p}^p\Big)^{1/p}$$ (actually if one is interested only in quadratic forms, it is enough to prove it for polynomials of degree 2). We do not know if such inequality is satisfied under the assumptions of Theorem \ref{thm:polynomial-Ising}.
\end{remark}

\begin{example}
Let us now present an example of a matrix $A$ for which Theorem \ref{thm:quadratic-forms-improved} gives a substantially better tail estimate that the one given by the Hanson-Wright inequality. One possibility is to tensorize
Example \ref{ex:easy-sum}, i.e., to consider the matrix $A_n = (a_{ij})_{i,j=1}^n$ given by $a_{ij} = \frac{1}{ij}$ for large values of $n$. Noting that $\langle AX,X\rangle = \langle v,X\rangle^2$ for $v = (1,1/2,\ldots,1/n)$ one can argue that this example is still rather about linear combinations than quadratic forms. Let us therefore leave the details to the Reader and instead consider the matrix $A_n$ given by $a_{ij} = \frac{1}{(i+j)^2}$. It is easy to see that for any $n$, $A_n$ is positive definite (e.g., by noting that for a standard exponential variable $Y$, and $t \ge 0$ we have $\E e^{-tY} = \frac{1}{1+t}$ and using basic properties of the Laplace transform). Now both $\|A_n\|_{HS}$ and $\|A_n\|_{\ell_2^n\to \ell_2^n}$ are of order $\Omega(1)$ as $n \to \infty$, and so, if $X_n$ is a sequence of centered random vectors in $[-1,1]^n$ satisfying \eqref{eq:convex-concentration-mean} with $K$ independent of $n$, then the Hanson-Wright type inequalities \eqref{eq:Hanson-Wright-123} give
\begin{displaymath}
  \p(|\langle A_n X_n,X_n\rangle - \E \langle A_nX_n,X_n\rangle| \ge t) \le 2\exp(-ct)
\end{displaymath}
for some dimension independent constant $c$. On the other hand, it is not difficult to check that we have
$\|A_n\|_{\{1,2\},p} \le C$ and $\|A_n\|_{\{1\}\{2\},p} \le C\log p$
for some dimension-independent constant $C$.
Thus Theorem \ref{thm:quadratic-forms-improved} gives
\begin{equation}\label{eq:final-example}
\p(|\langle A_n X_n,X_n\rangle - \E \langle AX_n,X_n\rangle| \ge t) \le 2e^{-e^{c't}},
\end{equation}
where $c'$ is another dimension independent constant and we again obtain a strengthened Gumbel type behavior.

The above example is primarily an illustration of the difference between the norms $\sqrt{p}\|A\|_{HS} + p\|A\|_{\ell_2^n\to \ell_2^n}$, used in the Hanson-Wright inequality and the norms $\|A\|_{\{1,2\},p} + \|A\|_{\{1\}\{2\},p}$ used in estimates of Rademacher type given in Theorem \ref{thm:quadratic-forms-improved}, but in fact one can recover \eqref{eq:final-example} by splitting appropriately the matrix $A$ into a sum of two matrices, applying the Hanson-Wright inequality to one of them and the trivial pointwise bound to the other one (similarly as in the one-dimensional case of Example \ref{ex:easy-sum}, where one can apply the pointwise bounds together with the Khintchine inequality). This strategy is however limited, as in general there does not exist a constant $C$, independent of $p$ and $n$ such that for all $n\times n$ matrices $A$ and $p\ge 2$,
\begin{displaymath}
  \inf\{\sqrt{p}\|B\|_{HS} + p\|B\|_{\ell_2^n\to \ell_2^n} + \|D\|_{\ell_1^n(\ell_1^n)} \colon B+D = A\} \le C(\|A\|_{\{1,2\},p} + \|A\|_{\{1\}\{2\},p}),
\end{displaymath}
where $\|(d_{ij})_{ij}\|_{\ell_1^n(\ell_1^n)} = \sum_{i,j\le n} |d_{ij}|$. To see this one can consider e.g., the matrix given in Remark \ref{rem:decomposition-fails} or matrices of the form $A = vv^T$, where $v$ has one coordinate equal to $1$ and the remaining ones equal to $p/n$, and $p \to \infty$ with $n$ at an appropriate speed (we leave the details to the Reader). This shows that estimates of the form \eqref{eq:upper-estimates-tails} do improve on the Hanson-Wright inequality.
\end{example}

\begin{proof}[Proof of Theorem \ref{thm:quadratic-forms-improved}]
Denote $f(x) = \langle Ax, x\rangle$. By our assumptions this is a convex function and therefore, similarly as in the proof of Proposition \ref{prop:vector-valued}, we can write for any $x,y \in [-1,1]^n$,
\begin{displaymath}
  f(x) - f(y)  \le  \langle \nabla f(x), x - y\rangle \le (2 + \frac{1}{\sqrt{p}}|x-y|)\|\nabla f(x)\|_{\{1\},p}.
\end{displaymath}
Therefore (denoting $M = \Med f(X)$), if  $f(x) > M + 3\|\nabla f(x)\|_{\{1\},p}$ then
$$g(x) := {\rm dist}\Big(x,\{y\in[-1,1]^n\colon f(y) \le M\}\Big) > \sqrt{p}$$ and thus by the convex concentration assumption applied to the function $g$ (note that $g$ is convex and $\Med g(X) = 0$), we get
\begin{displaymath}
  \p(f(X) - M \ge 3\|\nabla f(X)\|_{\{1\},p}) \le 2e^{-p/K^2}
\end{displaymath}
(observe that here we bound $f(X) - M$ by a random quantity).

As $\nabla f(X) = 2 A X$, we can apply Corollary \ref{cor:vector} to the norm $\|\cdot\|_{\{1\},p}$ to obtain
\begin{displaymath}
  \p(\|\nabla f(X)\|_{\{1\},p} \ge 2\E \|AX\|_{\{1\},p} + C \|A\|_{\{1\}\{2\},p} )\le 4e^{-p/K^2}.
\end{displaymath}
for some universal constant $C$.
Combining the two last inequalities we obtain for some (new) constant $C$,
\begin{align}\label{eq:quadratic-form-upper-tail}
  \p\Big(\langle AX,X\rangle  \ge M + C(\E\|AX\|_{\{1\},p} + \|A\|_{\{1\}\{2\},p} ) \Big)\le 6e^{-p/K^2}.
\end{align}

To get a bound on the upper tail (above the median) it thus suffices to estimate $\E \|AX\|_{\{1\},p}$. In \cite{MR1686370} Lata{\l}a proved that in the case when $Y$ is a vector of independent Rademacher variables,
\begin{displaymath}
  \E \|AY\|_{\{1\},p} \le C(\|A\|_{\{1,2\},p} + \|A\|_{\{1\}\{2\},p}).
\end{displaymath}
In \cite{MR3052405} it is mentioned (see the remark before Lemma 8.4) that this inequality can be proved by a chaining argument, relying on concentration properties of vector valued linear combinations of Rademacher variables. Since  $X$ is a centered random vector, satisfying analogous concentration properties as $Y$, one could follow this approach to recover the above inequality for $X$. However the formulations and proofs in \cite{MR1686370,MR3052405} are given for general independent  variables with log-concave tails and translating the arguments of \cite{MR3052405}, even if straightforward, is quite tedious. Therefore, instead we will use a recent deep result of Bednorz and Lata{\l}a \cite{MR3245015} concerning suprema of Rademacher processes, which will allow to directly reduce estimates for $X$ to the case of random signs.
Their Theorem 1.1 (in a finite dimensional formulation suitable for our purposes) asserts that if $T \subseteq \R^n$ then there exists a decomposition $T = T_1+T_2$ such that
\begin{displaymath}
  \sup_{t\in T_1} \sum_{i=1}^n |t_i| \le C\E \sup_{t \in T}\sum_{i=1}^n \varepsilon_i t_i
\end{displaymath}
and
\begin{displaymath}
  \E \sup_{t \in T_2} \sum_{i=1}^n g_i t_i \le C\E \sup_{t \in T} \sum_{i=1}^n \varepsilon_i t_i,
\end{displaymath}
where $\varepsilon_i$, $g_i$ are sequences of i.i.d.\ resp.\ Rademacher and Gaussian variables, and $C$ is a universal constant.

Since our $X$ satisfies convex concentration property, it is in particular subgaussian with constant $K$ and by another deep result, Talagrand's Majorizing Measure Theorem (see \cite{MR3184689}), we have
\begin{displaymath}
  \E \sup_{t \in T_2} \sum_{i=1}^n X_i t_i \le CK\E \sup_{t \in T_2} \sum_{i=1}^n g_i t_i
\end{displaymath}
for any set $T_2$. Therefore, expressing $\|AX\|_{\{1\},p}$ as a supremum of linear combinations of $X_i$'s and using the above estimates together with the inequality $|X_i| \le 1$, we see that
\begin{align}\label{eq:bound-on-expected-norm}
  \E \|AX\|_{\{1\},p} \le C(1+K)\E \|AY\|_{\{1\},p} \le C'(1+K)(\|A\|_{\{1,2\},p} + \|A\|_{\{1\}\{2\},p})
\end{align}
for some universal constants $C,C'$. Going back to \eqref{eq:quadratic-form-upper-tail} we obtain
\begin{align}
\p\Big(\langle AX,X\rangle  \ge M + C(1+K)(\|A\|_{\{1,2\},p} + \|A\|_{\{1\}\{2\},p} ) \Big)\le 6e^{-p/K^2}.
\end{align}
for some (new) universal constant $C$. Using the fact that $\|A\|_{\{1,2\},tp} \le C'' \sqrt{t}\|A\|_{\{1,2\},p} $ and $\|A\|_{\{1\}\{2\},tp} \le t \|A\|_{\{1\}\{2\},p}$ for $t \ge 1$ and some universal constant $C''$, we can easily replace the right hand side by $4e^{-p}$ at the cost of changing $C(1+K)$ to some constant $C_K$ (which can be clearly expressed explicitly in terms of $C$ and $K$).

As for the lower tail, \cite[Theorem 6.5 (ii)]{2017arXiv170301765A} implies that
under \eqref{eq:improved-convex-concentration}, for any smooth convex function $f\colon [-1,1]^n \to \R$,
\begin{displaymath}
  \p(f(X) \le \Med  f(X) - C_K \E\|\nabla f(X)\|_{\{1\},p}) \le 4e^{-p}
\end{displaymath}
(formally this theorem is stated for norms of the form $\|x\| = \sup\{\langle x,y\rangle\colon \theta(y) \le p\}$ with $\theta \colon \R^n \to \R$, whereas by \eqref{eq:dual-representation} $\|\cdot\|_{\{1\},p}$ corresponds to $\theta (y) = |y|^2$ if $|y|_\infty \le 1$ and $\theta(y) = \infty$ otherwise, but the proof given in \cite{2017arXiv170301765A} does not use the finiteness of $\theta$ and in fact one can also deduce the result by approximating $\theta$ with finite-valued functions).
In particular, the above inequality for $f(x) = \langle Ax,x \rangle$ gives
\begin{displaymath}
  \p\Big(\langle AX,X\rangle  \le M -C_K\E \|AX\|_{\{1\},p} \Big)\le 4e^{-p}.
\end{displaymath}

The inequality \eqref{eq:lower-estimates-tails} with the mean replaced by the median follows now by \eqref{eq:bound-on-expected-norm} and from the observation that
\begin{displaymath}
  \E\|AZ\|_{\{1\},p} \le \sqrt{p}\E|A X| \le \sqrt{p}\sqrt{\E |AX|^2} \le CK\sqrt{p}\|A\|_{\{1,2\}},
\end{displaymath}
where the last inequality follows easily by \eqref{eq:convex-concentration-mean} and integration by parts.

Now, \eqref{eq:upper-estimates-tails} and \eqref{eq:lower-estimates-tails} with the median instead of the mean yield
\begin{displaymath}
  \p(|\langle AX,X\rangle - M| \ge C_K t\|A\|_{\{1,2\}}) \le 4e^{-t},
\end{displaymath}
which by another integration by parts gives $|\langle AX,X\rangle - M\| \le C_K' \|A\|_{\{1,2\}}$. Since for $p \ge 2$, $\|A\|_{\{1,2\}} \le C\|A\|_{\{1,2\},p}$, this easily allows to pass from concentration around median to concentration around mean (at the cost of increasing the values of the constant $C_K$).
\end{proof}

\bibliographystyle{amsplain}	
\bibliography{IsingPolynomials}
\end{document}